\documentstyle[12pt]{article} 

\newcommand\qed{\hfill $\Box$ \vspace{5mm}}

\newcommand\pthree{{\bf {\rm P\/}^3}}
\newcommand\pfour{{\bf {\rm P\/}^4}}

\newcommand\lm{\lambda}
\newcommand{\binom}[2] 
{ \left(\begin{array}{c}#1\\#2\end{array}\right)}

\begin{document}

\title{A smooth surface in $\pfour$ not of general type has degree at most 46}
\author{Michele Cook}
\date{}
\maketitle

This is a continuation of the papers of Braun and Fl{\o}ystad [BF], 
Cook [C1] and Braun and Cook [BC] bounding the degree of smooth surfaces 
not of general type in $\pfour$ using Generic Initial Ideal Theory.  
We will use and expand on the results of [BC]
using Liaison Theory. We will  further restrict the configurations 
of the generic initial ideal of a generic hyperplane section of the 
surface by considering the geometric implications of having between one 
and five  generators of the generic initial ideal in degree much higher 
than the other generators. We prove the following

\vspace{3mm}
{\bf Theorem 1.}

{\em Let $S$ be a smooth surface of degree $d$ in $\pfour$ not of 
general type. Then $d \leq 46$. }

\section{Definitions and Introduction}

First  let us recall some of the language of generic ideal theory.  
(For a full treatment refer to [B] and [G].)
As we will mainly concern ourselves with a generic hyperplane 
section of a surface in $\pfour$ we will restrict our definitions to 
curves in $\pthree$. 

\vspace{2mm}
Let ${\rm \bf C}[x_0, x_1, x_2, x_3]$ be the ring of polynomials 
of $\pthree$ with the reverse lexicographical ordering $\succ$.
Let $C$ be a curve of degree $d$ in $\pthree$, with ideal $I_C$.  
After a generic change of basis, the monomial ideal of initial 
terms of elements of $I_C$ under $\succ$ is called the {\it generic 
initial ideal} of $C$ and denoted  ${\rm gin}(I_C)$. 

\vspace{2mm}
Let $\Gamma$ be a generic hyperplane section of $C$. Then
the generic initial ideal of $\Gamma$, 
${\rm gin}(I_{\Gamma})$,  is an ideal in ${\rm \bf C}[x_0, x_1, x_2]$
and can be defined by 
$$ {\rm gin}(I_{\Gamma}) = ({\rm gin}(I_C)|_{x_3 =0})^{{\rm sat}}$$

where the saturation is with respect to $x_2$. 

As $I_C$ and $I_{\Gamma}$ are saturated ideals, 
${\rm gin}(I_C)$ is generated 
by monomials of the form $x_0^{i}x_1^{j}x_2^{k}$
and the  generic initial ideal of $\Gamma$ is of the form 
$$ {\rm gin}(I_{\Gamma}) = (x_0^s,\  x_0^{s-1}x_1^{\lambda_{s-1}},
 \ \dots , \  x_1^{\lambda_0}), $$
where $\sum \lambda_i = d$ and $\lambda_i \geq \lambda_{i+1}+1$. 
Furthermore, as the points of $\Gamma$ are in uniform 
position, the work of Gruson and Peskine [GP] tells us that 
$\lambda_{i+1}+2 \geq \lambda_i $. 
The $\{ \lm_0, \  \lm_1, \  \dots   , \ \lm_{s-1} \}$
 are called the {\it connected invariants}  of $\Gamma$ or $C$.

\vspace{3mm}
{\bf Definition}.
A monomial $x_0^ax_1^bx_2^c$ is a {\it sporadic zero} of $C$ 
if $x_0^ax_1^bx_2^c \notin {\rm gin}(I_C)$, 
but $x_0^ax_1^b \in {\rm gin}(I_{\Gamma}$). 
(i.e. there exists $c' > c$ such that 
$x_0^ax_1^bx_2^{c'} \in {\rm gin}(I_C)$. )

\vspace{2mm}
Notice that every generator $x_0^ax_1^bx_2^{c'}$ of ${\rm gin}(I_C)$
with $c' > 0$ gives rise to sporadic zeros $x_0^ax_1^bx_2^c$ for all
$0 \leq c < c'$.

\vspace{2mm}
In [BF], it was shown that if $S$ is a smooth surface of degree $d$ 
in $\pfour$ not of general type, whose generic hyperplane section, 
$C$, has $\alpha_t$  sporadic zeros  in degree $t$ and  invariants 
 $\{ \lm_0, \  \lm_1, \  \dots   , \ \lm_{s-1} \}$,  then the following 
equation must be satisfied

\begin{eqnarray}
 0 \geq & d^2-&5d-18-10\sum_{i=0}^{s-1}(\binom{\lm_i}{2}+
(i-1)\lm_i) + \nonumber \\
 & & 12 \sum_{i=0}^{s-1} (\binom{\lm_i+i-1}{3} - 
\binom{i-1}{3})  - \sum_{t } \alpha_t (12t-22).
\end{eqnarray} 

\vspace{2mm}
For each degree $d$, there are only a few possibilities for  
$\{ \lm_0, \  \lm_1, \  \dots   , \ \lm_{s-1} \}$ and for each  
$d$ and $\{ \lm_0, \  \lm_1, \  \dots   , \ \lm_{s-1} \}$ there is 
an upper bound on the number of sporadic zeros by the work of Ellingsrud 
and Peskine [EP].
Thus, in order to bound $d$,  we need to find the smallest  
possible upper bound for  $A = \sum_{t=0}^{m} \alpha_t t$, for each 
$d$ and $\{ \lm_0, \  \lm_1, \  \dots   , \ \lm_{s-1} \}$. 
 
\vspace{2mm}
We will bound $A$ by bounding the 
degree of the generators of ${\rm gin}(I_C)$. We will show that 
if there are $a$ generators in a high degree $\geq  r > \frac{d}{2}$, 
and all the others are in degree $\leq r-2$, then $a \geq 6$.
Thus we can find an upper bound for $A$ by assuming that either all 
generators of ${\rm gin}(I_C)$ are in degree $\leq \frac{d}{2}$  
or if there exists a generator in degree 
$r > \frac{d}{2}$ then either there is a generator in degree $r-1$
or there are six generators in degree $\geq r$. Having done this we 
find that $d \leq 48$.  We will then use the connectedness of the 
invariants of $C$ (see [C2]) to lower  the bound to $46$ by considering 
the last few examples on a case by case basis. We will also give an example 
of a connected Borel-fixed monomial ideal in degree $46$ where 
all the conditions of this paper are satisfied. (However, whether or not 
this example corresponds to a actual curve is another matter.)

\section{Restricting the generators in high degree}

Let $S$ be a smooth surface not of general type in $\pfour$ with 
generic hyperplane section $C$. Let 
$\{ \lm_0, \  \lm_1, \  \dots \ , \ \lm_{s-1} \}$ be the invariants 
of $C$.  Assume that $d > (s-1)^2 +1$. In [BC] we found that 
$s \leq 7$ and $d \leq 66$. Furthermore,  if $s \leq 3, d \leq 8$ 
and  if $s=6$ or $7$,  $d \leq 44$. 
Thus we may restrict ourselves to the case $s = 4$ or $5$.

Let  $\gamma $ be  the number 
of sporadic zeros of $C$.  [EP] show that 
 
\begin{eqnarray}
{\rm for} \ s=4 & \gamma \leq 
1 + \sum_{i=0}^{s-1}(\binom{\lm_i}{2}+(i-1)\lm_i) 
-\frac{d^2}{8} + \frac{9d}{8}, \\ 
{\rm for} \ s=5 &  \gamma \leq  1 + 
\sum_{i=0}^{s-1}(\binom{\lm_i}{2}+(i-1)\lm_i) 
-\frac{d^2-5d+10}{10}. 
\end{eqnarray}  

\vspace{2mm}
For equation (1) to hold for large degree, $A$ will need to 
be  large. Every generator of ${\rm gin}(I_C)$ of the form 
$x_0^ax_1^bx_2^c$ with $c > 0$ gives a sporadic zero in each 
degree $n$ for $a + b \leq n \leq a+ b+ c -1$. 
So, one would like generators of ${\rm gin}(I_C)$
in as high a degree as possible. 

We saw in [BC], that if there were only one generator in degree 
$r > \frac{d}{2}$ and all others are in degree $\leq r-2$, this 
would imply that there were an secant line of $C$ of order $r$ 
which leads to a contradiction if $d > 50$. In Appendix A we give a 
slight improvement on this Lemma, showing that one can obtain a 
contradiction for $d > 42$. (Although 42 is  not optimal, in that we 
could lower this bound on $d$ if we were more careful, it is quite 
sufficient for our needs.) Thus, if there is one generator in degree 
$r > \frac{d}{2}$,  there must be another in degree $\geq r-1$. 
We will continue with this argument and show that if there 
are $a$ generators in degree $\geq r > \frac{d}{2}$ and 
all the others are in degree $\leq r-2$, then $a \geq 6$.

\vspace{2mm}
Consider the following situation. Let $C$ be a 
non-degenerate curve in $\pthree$  of degree $d$. 
Suppose ${\rm gin}(I_C)$ has $a$ generators 
in  degree $\geq r > \frac{d}{2}$, where  $1 \leq a \leq 5$
and the rest  in degree $\leq r-2$. (Assume from this point, 
that  a generic 
change of coordinates has been made so that  
${\rm gin}(I_C) = {\rm in}(I_C)$,  the initial ideal of $C$.)

Let $J$ be the ideal generated by elements of $I_C$ in degree 
$\leq r-1$. Then the generators of ${\rm gin}(J)$ are the 
generators of ${\rm gin}(I_C)$ in degree $\leq r-2$. 
By considering the Hilbert function associated to $J$, 
one finds  that degree$(V(J)) = $ degree$\ (C)+a$. 
Hence, $V(J) = C \cup X$ and degree$(X) = a$.
Thus $X$ contains a pure one-dimensional scheme $Y$ of degree $a$. 

We will  show that $Y$ must either be non-reduced or reducible. 
Then given that $a \leq 5$, $Y$ must contain (perhaps with a 
multiple structure) a line or a conic. Furthermore this line 
(respectively conic) must meet $C$ in $> \frac{d}{2}$ 
(respectively $> d$) points (up to multiplicity). 
We will then proceed as in [BC] to show in Propositon 4 that 
$Y$ cannot contain such a line or such a conic. 

\vspace{3mm}
{\bf Proposition 2}

If $I_C$ has two generators of degree $m$ and $n$ so 
that $m + n -2 \leq \frac{d}{2}$, then $Y$ is either non-reduced or 
reducible.

\vspace{2mm}
(The condition of Proposition 2 is easily satisfied, if our curves 
arise as generic hyperplane sections of surfaces not of  general type. 
Each of these curves lies on a surface of degree $s \leq 7$
(and in the cases we are interested in $s = 4, \ {\rm or} \  5$). 
By connectedness $\lambda_{s-1} \leq \frac{d}{s}-\frac{s-1}{2} < \frac{d}{s}$. Furthermore, as there are only a few {\it chains} of sporadic zeros, the curve will also lie on a surface of degree 
$\lambda_{s-1} +(s-1) +\epsilon$ where $\epsilon$ is small and 
hence $s + \lambda_{s-1} +(s-1) +\epsilon$ will be much less than 
$\frac{d}{2}$. )

\vspace{2mm}
{\bf Proof.}

Suppose, for a contradiction, that $Y$ is reduced and irreducible. 

\vspace{2mm}
{\bf Step 1}

$H^1({\cal I}_Y(d)) = 0$ for all $d \leq  0$.

\vspace{2mm}
Let $H$ be a generic hyperplane section of $\pthree$, then 
$Y \cap H$ is a set of points, $\Gamma' $, and we have the following 
short exact sequence

$$ O \rightarrow {\cal I}_Y(d-1)  \rightarrow 
 {\cal I}_Y(d) \rightarrow  {\cal I}_{\Gamma'}(d) \rightarrow O. $$

$H^0({\cal I}_{\Gamma'}(d)) = 0 $ for $d \leq 0$, therefore

$$ O \rightarrow H^1({\cal I}_Y(d-1))  \rightarrow 
 H^1({\cal I}_Y(d)) \rightarrow   $$

for $d \leq 0$ and hence
$h^1({\cal I}_Y(d-1))  \leq h^1({\cal I}_Y(d)) $ for $d \leq 0$.

Now, the restriction sequence 
$$ O \rightarrow {\cal I}_Y \rightarrow 
 {\cal O}_{\pthree} \rightarrow  {\cal O}_{Y}
 \rightarrow O $$

gives the exact sequence of  cohomology

$$ O \rightarrow H^0({\cal I}_Y)  \rightarrow 
 H^0({\cal O}_{\pthree}) \rightarrow  H^0({\cal O}_{Y})
 \rightarrow H^1({\cal I}_Y)  \rightarrow O. $$

As $Y$ is reduced and irreducible
$h^0({\cal O}_Y) = 1$. $ H^0({\cal I}_Y) =0 $ and hence 
$ H^0({\cal O}_{\pthree}) \cong  H^0({\cal O}_{Y})$. 
This implies that  $H^1({\cal I}_Y) = 0$ 
and hence $H^1({\cal I}_Y(d)) = 0$ for all $d \leq 0$. 

\vspace{2mm}
{\bf Step 2}

All the  generators of  
${\rm gin}(I_{C \cup X}) = {\rm gin}(J)$ are in degree 
$\leq t_0 \leq r-2$. We will show that all the  generators of  
${\rm gin}(I_{C \cup Y})$ are also in degree $\leq t_0 \leq r-2$.

(Recall that the maximum degree of the minimal generators of 
${\rm gin}(I_C)$ is the same as the regularity of $C$. (See [B]))

\vspace{2mm}
Let $H$ be a generic hyperplane. Then 
$(C \cup X) \cap H = (C \cup Y) \cap H = \Gamma$
and we have the following commutative diagram 

$$ \begin{array}{ccccccc} 
& 0 && 0 && 0 & \\ 
& \downarrow && \downarrow  && 
\downarrow & \\
O \rightarrow  & {\cal I}_{C \cup X}(t-1) & \rightarrow &
 {\cal I}_{C \cup X}(t) & \rightarrow &  {\cal I}_{\Gamma}(t) & 
\rightarrow O \\ 
 & \downarrow \alpha_1&& \downarrow \alpha_2 && 
\downarrow \alpha_3 & \\
O \rightarrow  & {\cal I}_{C \cup Y}(t-1) & \rightarrow &
 {\cal I}_{C \cup Y}(t) & \rightarrow &  {\cal I}_{\Gamma}(t) & 
\rightarrow O, 
\end{array}
$$
where $\alpha_3 $ is an isomorphism and $\alpha_1 $ and 
$\alpha_2$ are injections. All generators of ${\rm gin}(I_{C \cup X})$
are in degree $\leq t_0$, hence ${\cal I }_{C \cup X}$ is $t_0$-regular
and $H^1({\cal I}_{C \cup X}(t)) = 0$ for $t \geq t_0-1$. 

Taking cohomology of the diagram above we get for 
$t \geq t_0-1$

$$ \begin{array}{ccccccc} 
& 0 & & 0 & \\ 
& \downarrow  & & \downarrow & & &\\
\rightarrow &
 H^0({\cal I}_{C \cup X}(t+1)) & \stackrel{\alpha}{\rightarrow} & 
H^0({\cal I}_{\Gamma}(t+1)) & 
\rightarrow & O & \\ 
 & \downarrow \gamma && 
\downarrow \delta & \\
\rightarrow &
 H^0({\cal I}_{C \cup Y}(t+1)) & \stackrel{\beta}{\rightarrow} &  
H^0({\cal I}_{\Gamma}(t+1))  & 
\rightarrow & H^1({\cal I}_{C \cup Y}(t)) & \rightarrow
\end{array}
$$

$\alpha$ is onto and $\delta$ is an isomorphism, therefore 
$\delta \circ \alpha$ is onto.  
$\beta \circ \gamma = \delta \circ \alpha$,  and 
so $\beta$ is onto. Hence the coker($\beta) = 0$ for all 
$t \geq t_0-1$. This means there are no sporadic zeros of 
$C \cup Y$ in degree $t \geq t_0$ and hence all generators 
of ${\rm gin}(I_{C \cup Y})$ are in degree $\leq t_0$. 

\vspace{2mm}
{\bf Step 3}

$C \subset c.i.(n, m) $ a complete intersection of type $(m, n)$, 
where $m$ and $n \leq r-2.$ Thus $C$ is linked via this complete 
intersection to a curve containing the curve $Y$. Let's call 
this curve $C' \cup Y$. 
$C'$ is also a curve (perhaps non-reduced or reducible) and it 
is linked via the same complete 
intersection to $C \cup Y$. 

Now 
$$    O \rightarrow   {\cal O}_{C' \cup Y} \rightarrow 
 {\cal O}_{C'} \oplus {\cal O}_{Y}  \rightarrow   {\cal O}_{C' \cap Y}
\rightarrow O 
$$
and hence for all $t$ 
 $$   H^0({\cal O}_{C' \cup Y}(t)) \hookrightarrow 
 H^0({\cal O}_{C'}(t)) \oplus H^0({\cal O}_{Y}(t)).  $$

Also for  an arbitrary scheme $S \subset \pthree$ 

$$    O \rightarrow   {\cal I}_{S} \rightarrow 
 {\cal O}_{\pthree}  \rightarrow   {\cal O}_{S}
\rightarrow O 
$$

and hence for $t < 0$, $H^0({\cal O}_S(t)) \cong H^1({\cal I}_S(t))$

By Step1, $H^1({\cal I}_Y(t)) = 0$ for $t\leq 0$ and hence 

 $$   H^0({\cal O}_{C' \cup Y}(t)) \hookrightarrow 
 H^0({\cal O}_{C'}(t))  \ { \rm for \ all } \  t < 0 $$

or equivalently 
 \begin{eqnarray}
   H^1({\cal I}_{C' \cup Y}(t)) \hookrightarrow 
 H^1({\cal I}_{C'}(t))  \ { \rm for \ all } \  t < 0.
\end{eqnarray}

Now ${\rm gin}(I_C)$ has a generator in degree $r$ and hence 
$C$ has a sporadic zero in degree $r-1$ which means that 
$H^1({\cal I}_C(r-2))\neq 0$. $C$ is linked to $C' \cup Y$ 
via a complete intersection of type $(m,n)$ 
and so by a fundamental theorem from liaison theory ([PS]) 
$$H^1({\cal I}_C(t)) \cong H^1({\cal I}_{C' \cup Y}(m + n - 4 - t)){\check{}}. $$
Hence 
$H^1({\cal I}_{C' \cup Y}(m + n - 4 - (r-2))) \neq 0$. 

On the other hand $C'$ is linked to $C \cup Y$ via a 
complete intersection of type $(m,n)$ and all the sporadic zeros of  
$C \cup Y$ are in degree $\leq t_0 -1 \leq r-3$. Hence 
$H^1({\cal I}_{C \cup Y}(t)) = 0 $ for all $t \geq t_0-1$. 
Again by the relationship between the deficiency modules 
of linked curves 
$$H^1({\cal I}_{C \cup Y}(t)) \cong H^1({\cal I}_{C'}(m + n - 4 - t)){\check{}}$$

and so $H^1({\cal I}_{C'}(m + n - 4 - t))= 0$ for all $t \geq t_0-1$.
$r-2 \geq t_0-1$, therefore $H^1({\cal I}_{C'}(m + n - 4 - (r-2)))= 0$.
Furthermore, by the assumptions of the proposition  $m+n-r-2 < 0$.  
This contradicts the injection in equation (4). 
\qed 

Let $I_C = (J, f_1, \  \dots,  \  f_b)$ where $ 1 \leq b \leq a \leq 5$ 
and $\frac{d}{2} < r = $ deg$(f_1) \leq \ \dots \ \leq$ deg$(f_b)$. 
As $Y$ is either non-reduced or reducible and  the degree of $Y$ is at most 5, $Y$ must contain 
either a line,  $l$,  or a conic, $Q$ perhaps with a multiple structure. 

\vspace{3mm}
{\bf Lemma 3}

There exists $i, 1 \leq i \leq b$ such that $l \cap F_i$ 
(respectively $Q \cap F_i$) in 
${\rm deg}(f_i) > \frac{d}{2}$ (respectively $2{\rm deg}(f_i) > d$)
points (up to multiplicity). 
(Here $F_i = \{x \in \pthree | f_i(x) = 0 \}$.)

\vspace{2mm}
{\bf Proof}

As the line $l$ or conic $Q$, with perhaps a  multiple structure, 
is contained in $V(J)$, then the 
reduced scheme will certainly be contained in $V(J)$ as the 
non-reduced structure is really information about the embedding. 
Thus $l$ or $Q$ is contained in $V(J)$.

If $b=1$, $I_C = (J, f_1)$ and  $l$ (respectively $Q$) 
$\notin F_1$. By Bezout's theorem  $l \cap F_1$ (respectively 
$Q \cap F_1$) in ${\rm deg}(f_1) > \frac{d}{2}$ \ 
(respectively $2 {\rm deg}(f_1) > d$) points (up to multiplicity). 

If $b >1$, consider the ideals 
$I_j = (J, f_1, \dots, \ \widehat{f_j}, \dots ,  \ f_b )$, 
$j = 1, \dots  , b$. 
As ${\rm in}(f_i)$ is of the form $x_0^{i_0}x_1^{i_1}x_2^{i_2}$  
with $i_2 > 0$,  deg$V(I_j) > $ deg$(C)$. Hence $V(I_j)$ contains 
a pure one dimensional scheme $Y_j$.
$C \subset C \cup Y_j \subset C \cup Y$  and hence $Y_j$ 
contains some component of $Y$. 

If $Y_j$ contains
a line or a conic (with possible multiple structure) we are done 
as in the case of $b=1$, replacing $J$ with $I_j$. 

The only other possibility is that $Y_j$ only 
contains the unique cubic or quartic (possibly) contained in $Y$. 
However this is only possible for at most one 
choice of $j$, (otherwise the cubic or quartic is contained in $C$.)  
Therefore there exists some $j$ such that $Y_j$ contains 
a line or a conic.
\qed

\vspace{3mm}
{\bf Proposition 4}

Let $S$ be a smooth surface not of general type of degree $d$ 
in $\pfour$ with generic hyperplane section $C$.
Suppose ${\rm gin}(I_C)$ has at least four generators, at least 
two in degree $\leq r-2$ and $a$ in degree $\geq r > \frac{d}{2}$.
Then $d \leq 42$ or $a \geq 6$.

\vspace{2mm}
{\bf Proof}

The proof is in two parts. The first uses the (slightly improved) 
result of [BC], which can be found at the end of this paper, 
which states
that if $S$ is  a smooth surface not of general type of degree $d$ 
in $\pfour$ with generic hyperplane section $C$ and $d > 42$, 
then $C$ cannot have a secant line of order $r > \frac{d}{2}.$

\vspace{2mm}
The second deals with the conic.
Let $Q$ be the conic as above. Let 
deg$(f_i) = t \geq r >  \frac{d}{2}$ and $F = \{ f_i=0 \}$. 
We may assume $Q$ is reduced and irreducible 
and hence (as we are in $\pthree$) a smooth conic. 
By Bezout's Theorem $ F \cap Q $ in $2t >d$ points 
(up to multiplicity) and all these points must lie on $C$. Let 
$F \cap Q = \sum m_ip_i $ where $p_i$ are the points of $C$. 

\vspace{2mm}
{\bf Claim}

$Q$ meets $C$ at $p_i$ with multiplicity $m_i$. 

\vspace{2mm}
Given the claim, as $Q$ is a conic $Q \subset P$ a plane 
in $\pthree$ and so $C \cap P$ in at least $2t > d$ points. But 
this means that $C$ is contained in the plane which is not possible.

\vspace{2mm}
{\bf Proof of the Claim}

The definition of intersection multiplicity given in [H] pg 427, is
as follows.
Let $X$ and $Y$ be varieties (in $\pthree$ for simplicity) 
meeting at a point $p$,  
then the intersection multiplicity 
$$i(X, Y, p) = \sum_{j} (-1)^j 
{\rm length}({\rm Tor}_j^A(A/{\bf a}, A/{\bf b})). $$

where $A ={\cal O}_{p, \pthree}$ is the local ring of $p$ and 
${\bf a}= I_X$ and ${\bf b} = I_Y$ are the ideals of $X$ and $Y$ 
in $A$. 

As $C$ is smooth, it is a local complete intersection, locally cut out 
at $p$ by coprime polynomials $g_1$ and $g_2$. We may assume
that these polynomials correspond to polynomials in $I_C$ of 
degree $t$.  
Suppose the conic $Q$ meets $V(g_1)$ and 
$V(g_2)$ at $p$ with multiplicity $m$ and thus 
$$ m = \sum_{j} (-1)^j {\rm length}({\rm Tor}_j^A(A/{\bf a}, A/{\bf b}_i))$$

Where ${\bf a} = I_Q$ and ${\bf b}_i = (g_i)$ for $ i = 1, 2$. 

We have the following short exact sequence

$$ 0 \rightarrow \frac{{\bf b_2}}{{\bf b_1}} \rightarrow 
\frac{A}{{\bf b_1}} \rightarrow \frac{A}{(g_1, g_2)} 
\rightarrow 0$$ 

and hence via the long exact sequence of Tor, 
$$\begin{array}{ll}
m & = \sum_j (-1)^j {\rm length}({\rm Tor}_j^A(\frac{A}{{\bf a}})), 
\frac{A}{{\bf b_1}}) \\
 & = \sum_j (-1)^j{\rm length}({\rm Tor}_j^A(\frac{A}{{\bf a}}, \frac{A}{(g_1,g_2)}))
+\sum_j (-1)^j {\rm length}({\rm Tor}_j^A(\frac{A}{{\bf a}}, 
\frac{{\bf b_2}}{{\bf b_1}})). 
\end{array}$$

Thus we need to show that 
$$\sum_j (-1)^j {\rm length}({\rm Tor}_j^A(\frac{A}{{\bf a}}, 
\frac{{\bf b_2}}{{\bf b_1}}) )= 0. $$

Now, 
$$ O \rightarrow {\bf b_1} \cap {\bf b_2} \rightarrow 
{\bf b_2} \rightarrow \frac{{\bf b_2}}{{\bf b_1}} 
\rightarrow O $$

Thus we need to show that 
$$\sum_j (-1)^j {\rm length}({\rm Tor}_j^A(\frac{A}{{\bf a}}, 
{\bf b_1} \cap {\bf b_2}) )= \sum_j (-1)^j {\rm length}({\rm Tor}_j^A(\frac{A}{{\bf a}}, {\bf b_2}) ). $$

${\bf b_1} \cap {\bf b_2} = (g_1g_2)$ and 
${\bf b_1} = (g_1)$ are principal ideals and hence have 
minimal resolutions 

$$ O \rightarrow A \rightarrow {\bf b_1} \cap {\bf b_2} 
\rightarrow O $$

$$ O \rightarrow A \rightarrow {\bf b_1} 
\rightarrow O $$

Therefore for any $A$-module $B$,  
${\rm Tor}^A_i(B, {\bf b_1} \cap {\bf b_2}) = 
{ \rm Tor}^A_i(B, {\bf b_1} ) = 0$ if $i \geq 1$. 
Moreover 
${\rm Tor}^A_0(\frac{A}{{\bf a}}, {\bf b_1} \cap {\bf b_2})  
= \frac{A}{{\bf a}} \otimes {\bf b_1} \cap {\bf b_2}
= \frac{{\bf b_1} \cap {\bf b_2}}{{\bf a} \cap {\bf b_1} \cap {\bf b_2}}$ 
and 
 ${\rm Tor}^A_0(\frac{A}{{\bf a}}, {\bf b_1})  
= \frac{A}{{\bf a}} \otimes {\bf b_1}
= \frac{{\bf b_1}}{{\bf a} \cap {\bf b_1}}.$ 

Let $\Phi : \frac{{\bf b_1}}{{\bf a} \cap {\bf b_1}}
\rightarrow \frac{{\bf b_1} \cap {\bf b_2}}{{\bf a} \cap {\bf b}_1 \cap {\bf b_2}}$  be defined by 
$\Phi(g_1a + {\bf a}  \cap {\bf b_1}) = g_1g_2a + {\bf a} \cap {\bf b_1}\cap 
{\bf b_2}$. This is an isomorphism, hence we are done.

\vspace{2mm}
Thus if $a \leq 5$, $Y$ either contains a secant line of $C$ of order 
$t > \frac{d}{2}$ which contradicts [BC] if $d > 42$ or $Y$ contains a 
conic meeting $C$ in $>d$ points which again leads to a contradiction. 
Therefore either $d \leq 42$ or $a \geq 6$.
\qed

\vspace{2mm}
Note that if $Y$ were of degree 6, there is the possibility 
that $Y$ is a union of two cubics, or a double cubic, which cannot 
be eliminated by the above methods. However, given the bounds 
on the number of sporadic zeros it would be very difficult 
to have more than $6$ generators in degree $> \frac{d}{2}$. And as we will see in the last example this is not the obstruction to lowering the bound further. 

\section{Bounding A}

 \vspace{3mm}
Thus if the degree of $S > 42$, a generic hyperplane section, 
$C$, of $S$ cannot have an $t$-secant line with $t > \frac{d}{2}$ 
nor a conic meeting the curve in $2t > d$ points. 
In terms of the generic initial ideal of $C$, this means that 

either 

(1) all generators of ${\rm gin}(I_C)$ are in degree $\leq \frac{d}{2}$ 

or 

(2) if there exists a generator of ${\rm gin}(I_C)$  in degree 
$r > \frac{d}{2}$ then either there exists a generator in degree 
$r-1$ or there exist a five more  generators in degree $\geq r$.

\vspace{2mm}
The idea now is to build Borel-fixed monomial ideals 
with as many generators in degree 
$\geq \lfloor \frac{d}{2} \rfloor$ as possible.

If  there can be 6 generators in degree 
$\geq \lfloor \frac{d}{2} \rfloor+1$, let the ideal have 
5 generators in degree 
$\lfloor \frac{d}{2} \rfloor +1 $ and the sixth generator be in 
as high a degree as possible using all the remaining sporadic zeros. 
(This configuration is unlikely to arise from a curve. 
However, for the purposes of Bounding $A$ it will give a good 
enough upper bound.) If $s=5$ the bound
on the number of sporadic zeros means that it is impossible to have 
6 or more generators in degree $\lfloor \frac{d}{2} \rfloor +1$, and 
it is $s=5$ which is giving the upper bound on $d$. 

Otherwise let the ideal have one generator in degree  
$\lfloor \frac{d}{2} \rfloor$, one in degree 
$\lfloor \frac{d}{2} \rfloor +1 $ and so on until all the 
sporadic zeros have been used up. 

The rest of the calculations are done by computer. 
For each $d \leq 66$ (the bound we obtained in [BC]) we find all 
connected invariants $\{ \lambda_{0}, \dots , \lambda_{s-1} \}$ 
with $s = 4$ or $5$ using a program of Rich Liebling. 
We then use equations
 $(2)$ and $(3)$ to find the maximal number of sporadic zeros, $z$. 
We then use Mathematica to find the maximal $A$ using the criteria 
above and see which examples satisfy equation $(1)$.  
The following is a list of the examples with the highest degree 
the Mathematica program gives  which satisfy the criteria 
above and equation (1). 

\vspace{3mm}
${\bf  s = 4}$

For $s=4$  the only possibilities are in degree $\leq 46$. Those of 
highest degree being 

\vspace{2mm}
\begin{center}
\begin{tabular}{|| c | c | c | c | c | c || } \hline
\# &degree & $z$ & $\{ \lambda_i \} $ & A-bound & neg \\
\hline 
1 & 46 & 50 & 14, \ 12, \ 11, \ 9 & 921 & -2 \\
2 & 45 & 49 & 14, \ 12, \ 10, \ 9 & 882 & -18 \\
3 & 45 & 48 & 13, \ 12, \ 11, \ 9 & 854 & -2 \\
\hline
\end{tabular}
\end{center}

\vspace{3mm}
${\bf  s = 5}$

For $s=5$ the only possibilities are in degree $\leq 48$. Those of 
highest degree being 

\vspace{2mm}
\begin{center} 
\begin{tabular}{|| c | c | c | c | c | c || } \hline
\# & degree & $z$ & $\{ \lambda_i \} $ & A-bound & neg \\
\hline 
1 & 48 & 45 & 13, \ 11, \ 10, \ 8, \ 6 & 810 & -14 \\
2 & 47 & 44 & 13, \ 11, \ 9, \ 8, \ 6 & 770 & -20 \\
3 & 47 & 43 & 12, \ 11, \ 10, \ 8, \ 6 & 743 & -4 \\
4 & 46 & 44 & 13, \ 11,   \ 9, \ 7, \ 6 & 770 & -56 \\
 \hline
\end{tabular}
\end{center}

\vspace{3mm}
{\bf Notes:}

\vspace{1mm}
{\bf 1.} $z$ is the maximal number of sporadic zeros.

\vspace{1mm}
{\bf 2.} {\it A-bound} is the upper bound on $A$ found using the 
criteria above. 

\vspace{1mm}
{\bf 3.} The column {\it neg} gives the value of equation (1) divided 
by 12 then rounded up. Hence for each of the examples above in 
degree $\geq 47$ we 
need to show that $A < $A- bound $-$  neg in order to eliminate that 
particular case. 

\vspace{1mm}
{\bf 4.} The input list was not actually the full list of 
possible connected invariants for each degree. 
Certain configurations give rise to Arithmetically Cohen-Macaulay 
(ACM) curves, or equivalently, curves without sporadic zeros. 
If a  curve has no sporadic zeros then $A =0 $, then in order to 
satisfy equation (1) we must have 
for $s = 4$,  $d \leq 10$ and for $s= 5$, $d \leq 17$. 

Invariants which give ACM curves (for $s \geq 4$) are 

$(i)$ those whose consecutive invariants differ by two in which 
case $C$ is a complete intersection of type $(s, \lambda_{s-1}+s-1)$. 
For example the invariants $ \{ 6, 8, 10, 12, 14 \}$ correspond to 
a complete intersection of type $(5, 10)$.

\vspace{1mm}
$(ii)$ those whose  consecutive invariants differ by two 
except $\lambda_0 = \lambda_1 + 1$ in which case $C$ is linked to a 
line by a complete intersection of type $(s, \lambda_{s-1} + s-1)$. 
Then as a line is ACM, $C$ is ACM.  (See Rao [R] or 
Migliore [M]) For example the invariants $ \{ 6, 8, 10, 12, 13 \}$ 
correspond to a
curve linked by a complete intersection of type $(5, 10)$ to a line.

\vspace{1mm}
{\bf 5.} We will show that example 4 for $s=5$ gives rise to  
a Borel-fixed, connected monomial ideal
which satisfies all the conditions of this paper.

\vspace{2mm}
{\bf Eliminating} ${\bf (13, \ 11, \ 10, \ 8, \ 6\ )}$

\vspace{1mm}
We will eliminate the degree 48 configuration using 
connectedness. The two configurations in degree 47 can 
also be eliminated in a similar way.

\vspace{2mm}
If $x_0^4x_1^6$ is a generator of ${\rm gin}(I_C)$ 
then $C$ is contained in a complete intersection of type $(5, 10)$ and 

$$ {\rm gin}(I_C) \supseteq (x_0^5, x_0^4x_1^6, x_0^3x_1^8, x_0^2x_1^{10}, x_0x_1^{12}, x_1^{14}). $$

This means that the only monomials of the form 
 $x_0^ax_1^bx_2^c$ with $c >0$   which can be generators of 
${\rm gin}(I_C)$ are those with $(a, b) = (1, 11)$ or $(0, 13)$.
Then the best one could hope for is that one is of degree 24, 
and the other is of degree 25 and
$$ A \leq \sum_{t = 12}^{23} t + \sum_{t=13}^{24}t = 156 < 810-14.$$ 
So this is not a possibility. 

\vspace{1mm}
Therefore $x_0^4x_1^6$ is a sporadic zero. 
As $z = 45$, there are at most three 
generators in degree $\geq 24$ and taking connectedness 
into account (See [C2]) there are only a few possibilities. The one giving 
the best bound on $A$ would be if $x_1^{13}x_2^{12}$, 
$x_0^1x_1^{11}x_2^{12}$ and  $x_0^4x_1^6x_2^{16}$ are generators. 
this leaves 5 sporadic zeros which one can use up by making 
$x_0^3x_1^8x_2^{5}$ a generator. This gives an upper bound on $A$ 
of $777 < 810-14$, therefore this particular set of invariants is not 
possible. 

\vspace{2mm}
{\bf Example} ${\bf (13,\ 11, \ 9, \ 7, \ 6\ )}$

\vspace{1mm}
The degree 46 configuration for $s=5$ is actually possible. 
(i.e. one can create a Borel-fixed monomial ideal which is connected 
and still get an upper bound on $A$ which is big enough.)

Let ${\rm gin}(I_C)$ be defined as follows

$$\begin{array}{ll}
{\rm gin}(I_C) = (x_0^5, \  x_0^4x_1^6x_2,  &x_0^4x_1^7, \  
x_0^3x_1^7x_2^{13}, \ x_0^3x_1^8, \ x_0^2x_1^9x_2^{13}, \\ 
&x_0^2x_1^{10}, \ x_0x_1^{11}x_2^{13}, \ x_0x_1^{12}, \ 
x_1^{13}x_2^4, \ x_1^{14})
\end{array}$$

(Notice that $x_0^4x_1^6$ cannot be a generator of ${\rm gin}(I_C)$
otherwise $C$ would be linked to a plane quartic 
and hence Arithmetically Cohen-Macaulay.)

In this case we get $A = 731 > 770-56$.

\vspace{3mm}
{\it Acknowledgments}. 
This work was completed for the Summer School in Commutative 
algebra at the Centre de Recerca  Matem\`{a}tica, Barcelona in the 
Summer of 1996. I would like to thank the Centre for allowing me to 
participate. I would also like to thank the Association for
 Women in Mathematics for their generous contribution toward 
my travel expenses to the meeting. 

\section*{Appendix A}

We saw in [BC], that if $S$ is a smooth surface not of general type 
in $\pfour$, whose hyperplane section  $C_h$ has an secant line 
of order $r$, with $r > \frac{d}{2}$, then $S$ contains a plane 
curve of degree $r > \frac{d}{2}$. In [BC] we showed that this was 
impossible if $d > 50$, here we will improve the lemma slightly to 
show that 
this is impossible if $d >42$. 

\vspace{3mm}
{\bf Lemma}
{\em If $S$ is a smooth surface not of general type in $\pfour$ of 
degree $ d > 42$, $S$ cannot contain a plane curve of degree 
$r > \frac{d}{2}$. }

\vspace{2mm}
{\bf Proof}

We may assume that $d > (s-1)^2+1$ then using the bounds

\begin{eqnarray}
1+ \sum_{i=0}^{s-1}(\binom{\lm_i}{2}+(i-1)\lm_i) & \leq & 
\frac{d^2}{2s}+(s-4)\frac{d}{2}+1 
\end{eqnarray} 

\begin{eqnarray}
\sum_{t=0}^{s-1} (\binom{\lm_t+t-1}{3} - \binom{t-1}{3}) & \geq &
s\binom{\frac{d}{s}+\frac{s-3}{2}}{3}+1-\binom{s-1}{4}
\end{eqnarray} 

found in [GP] and [BF] (respectively),  equation (1) can be approximated by 

\begin{eqnarray}
  0 & \geq  & d^2-5d-18-10(\frac{d^2}{2s}+(s-4)\frac{d}{2}) 
+12 s\binom{\frac{d}{s}+\frac{s-3}{2}}{3} \nonumber \\
&& + 12(1-\binom{s-1}{4}) - \sum_{t =0}^{m} \alpha_t (12t-22).  
\end{eqnarray}

Let us first find a lower bound for the number of sporadic 
zeros of a generic hyperplane section of $S$. 

For $s=4 $, suppose that the number of sporadic zeros is 
$\leq \frac{3d}{4}$ then,  naively,  
$A \leq \sum_{\lambda_0}^{\lambda_0 + \frac{3d}{4} -1} t.$ 
By connectedness $\lambda_0 \leq \frac{d}{4}+3$ and hence 
$A \leq \frac{5}{32}d^2 + \frac{13}{8}d-3$. Substituting back into 
equation $(7)$ we get 
$$ 0 \geq \frac{d^3}{8}-\frac{23}{8}d^2-\frac{17}{2}d+33 $$
and hence $d \leq 25$. Therefore we may assume that the 
number of sporadic zeros is $ >  \frac{3d}{4}$, then 
$ g(C_h)  < \frac{d^2}{8}  +1 -\frac{3d}{4}$.

Similarly for $s=5 $, suppose that the number of sporadic zeros is 
$\leq \frac{2d}{5}$ then
$A \leq \sum_{\lambda_0}^{\lambda_0 + \frac{2d}{5} -1} t.$ 
By connectedness $\lambda_0 \leq \frac{d}{5}+4$ and hence 
$A \leq \frac{4}{25}d^2 + \frac{7}{5}d$. Substituting back into 
equation $(7)$ we get 
$$ 0 \geq \frac{d^3}{25}-\frac{24}{25}d^2-10d-9 $$
and hence $d \leq 35$. Therefore we may assume that the 
number of sporadic zeros is $ >  \frac{2d}{5}$, then 
$ g(C_h)  < \frac{d^2}{10} +\frac{d}{2} +1 -\frac{2d}{5}$.

\vspace{2mm}
Let $C \subset P$ be a plane curve  of degree 
$r > \frac{d}{2}$ contained in $S$. Let $H$ be a hyperplane containing $C$. 

Then
$S \cap H = C_h = C \cup C_{res} $. 

We have 
$$ 0 \rightarrow {\cal O}_{C \cup C_{res}} \rightarrow 
{\cal O}_{C } \oplus {\cal O}_{C_{res}} \rightarrow {\cal O}_{C \cap C_{res}} 
\rightarrow 0 $$

therefore 
$$h^1 ({\cal O}_{C_h}) \geq h^1({\cal O}_{C }) + 
h^1({\cal O}_{C _{res}})$$
and hence 
$$ g(C_h) \geq g(C) + g(C_{res}) \geq g(C). $$

$C$ is a plane curve of degree $d_C \geq \frac{d}{2}$ and so 
$$g(C) = \frac{(d_C-1)(d_C-2)}{2} - \delta 
\geq \frac{(\frac{d}{2}-1)(\frac{d}{2}-2)}{2} $$

Hence
$$ \begin{array}{ll}
{\rm for} \ s=4 &  \frac{d^2}{8}  +1 -\frac{3d}{4} >
 \frac{(\frac{d}{2}-1)(\frac{d}{2}-2)}{2}  \\ 
{\rm for} \ s = 5 &  \frac{d^2}{10} +\frac{d}{2} +1 -\frac{2d}{5} 
\geq \frac{(\frac{d}{2}-1)(\frac{d}{2}-2)}{2}   \\
{\rm for} \ s=6, 7 &  \frac{d^2}{2s} + (s-4)\frac{d}{2} +1 \geq
\frac{(\frac{d}{2}-1)(\frac{d}{2}-2)}{2}   
\end{array}
$$

This means for $s=4$ we have a contradiction and hence $d \leq 25$, 
for $s=5$, $d < 34$ and for $s =6$ and $7$, $d \leq 42$,
\qed

\section*{References} 

{\bf [B]} D. Bayer {\it The division algorithm and the Hilbert scheme}, 
Ph.D. Thesis, Harvard University (1982). 

{\bf [BC]} R. Braun, M. Cook {\it A smooth surface in $\pfour$
not of general type has degree at most 66}, to appear in Compositio
Mathematica.

{\bf [BF]} R. Braun, G. Fl{\o}ystad {\it A bound for the degree 
of smooth surfaces in $\pfour$ not of general type}, Compositio 
Mathematica, Vol. 93, No. 2, September(I) (1994) 211-229.

{\bf [BPV]} W. Barth, C. Peters, A. Van de Ven 
{\it Compact Complex Surfaces}, Springer-Verlag (1984). 

{\bf [C1]} M. Cook {\it An improved bound for the degree of smooth 
surfaces in $\pfour$ not of general type},  Compositio
Mathematica, Vol. 102, No. 2, June 1996,  141-145.

{\bf [C2]} M. Cook {\it The connectedness of space curve invariants}
to appear in Compositio Mathematica.

{\bf [EP]} G. Ellingsrud, C. Peskine {\it Sur les surfaces lisses de 
$\pfour$}, Invent. Math. 95 (1989) 1-11.

{\bf [G]} M. Green {\it Generic Initial Ideals}  Notes  from the Summer School on Commutative Algebra, Centre de Recerca Matam\`{a}tica, Spain, 1996. 

{\bf [GP]} L. Gruson, C. Peskine {\it Genres des courbes de l'espace 
projectif}, Lecture Notes in Mathematics, Algebraic Geometry, 
Troms{\o} 1977, 687 (1977) 31-59.

{\bf [H]} R. Hartshorne {\it Algebraic Geometry}, Springer-Verlag (1977). 

{\bf [K]} L. Koelblen {\it Surfaces de $\pfour$ trac{\'e}es sur une
 hypersurfaces cubique} Journal f{\"u}r die Riene and Angewandte 
Mathematik, 433 (1992) 113-141.

{\bf [M]} J. Migliore {\it An introduction to deficiency modules and 
liaison theory for subspaces of projective space} Seoul National University, Lecture Notes Series \# 24. 

{\bf [PS]} C. Peskine, L. Szpiro {\it Liaison des vari{\'e}t{\'e}s 
alg{|'e}brique I} Invent. Math. 26 (1974) 271-302.
 
{\bf [R]} P. Rao {\it Liaison among curves in $\pthree$.} Invent. 
Math. 50 (1979) 205-217.

\vspace{5mm}
Michele Cook

Department of Mathematics

Pomona College

610 N. College Avenue

Claremont, CA 91711-6348

e-mail mcook\verb+@+pomona.edu

\end{document}